\documentstyle{amsppt}
\voffset-10mm
\magnification1200
\pagewidth{130mm}
\pageheight{204mm}
\hfuzz=2.5pt\rightskip=0pt plus1pt
\binoppenalty=10000\relpenalty=10000\relax
\TagsOnRight
\loadbold
\nologo
\addto\tenpoint{\normalbaselineskip=1.\normalbaselineskip\normalbaselines}
\let\epsilon\varepsilon
\let\le\leqslant
\let\ge\geqslant
\let\wt\widetilde

\let\[\lfloor
\let\]\rfloor
\redefine\d{\roman d}
\define\Li{\operatorname{Li}}

\topmatter
\title
Approximations to -, di- and tri- logarithms
\endtitle
\author
Wadim Zudilin\footnotemark"$^\ddag$"\ {\rm(Moscow)}
\endauthor
\date
\hbox to100mm{\vbox{\hsize=100mm%
\centerline{E-print \tt math.NT/0409023}
\smallskip
\centerline{August 2004}
}}
\enddate
\address
Department of Mechanics and Mathematics,
Moscow Lomonosov State University,
Vorobiovy Gory, GSP-2,
119992 Moscow, RUSSIA
\endaddress
\email
{\tt wadim\@ips.ras.ru}
\endemail
\abstract
We propose hypergeometric constructions of simultaneous
approximations to polylogarithms. These approximations
suit for computing the values of polylogarithms and
satisfy 4-term Ap\'ery-like (polynomial) recursions.
\endabstract
\keywords
Dilogarithm, trilogarithm, zeta value, hypergeometric series,
algorithm of creative telescoping, polynomial recursion,
Ap\'ery's approximations
\endkeywords
\subjclass
33C20, 33F10, 11J70, 11M06
\endsubjclass
\endtopmatter
\leftheadtext{W.~Zudilin}
\footnotetext"$^\ddag$"{The work is partially supported
by grant no.~03-01-00359 of the Russian Foundation for Basic Research.}
\document

The series for the logarithm function
$$
\log(1-z)=-\sum_{n=1}^\infty\frac{z^n}n,
\qquad |z|<1,
$$
is so simple and nice that mathematicians immediately
generalize it by introducing the polylogarithms
$$
\Li_s(z)=\sum_{n=1}^\infty\frac{z^n}{n^s},
\qquad |z|<1, \quad s=1,2,\dots,
\tag1
$$
and considering then multiple, $q$-basic, $p$-adic
and any other possible generalizations, just to have
a serious mathematical research (i.e\. to have some fun).
It is so that nobody could now overview the whole
amount of the results around all these generalizations
of the logarithm, since the literature on the subject
increases to infinity as a geometric progression
(almost hypergeometrically).

It is not surprising
that the transcendence number theory also dreams
of getting new and new results for the values of
the polylogarithms, especially after Lindemann's
proof of the transcendence of $\log x$ for any
algebraic $x$ different from~$0$ and~$1$. The main
problems (or, if you like, intrigues) are therefore
extensions of the result to the values of~\thetag{1},
where Lindemann's argument based on proving the
transcendence for the inverse, exponential, function
does not work in an obvious manner any more.
And we now have, so far irrationality and linear
independence results for the polylogarithm values
at nonzero rational points close to zero,
thanks to contributions of Maier~\cite{Ma} (1927),
Galochkin~\cite{Ga} (1974), Nikishin~\cite{Ni} (1978),
Chudnovsky~\cite{Ch} (1979),
Hata~\cite{Ha1},~\cite{Ha2} (1990 and 1993),
Rhin and Viola~\cite{RV} (2004).
The very last piece of news is the irrationality
of $\Li_2(1/q)$ for $q$~integer, $q\le-5$ or $q\ge6$,
obtained by the powerful arithmetic method in~\cite{RV},
which improves the range of~\cite{Ha2} by adding $q=6$
(the work~\cite{RV} also includes quantitative improvements
of the irrationality in other cases, but we do not touch
this subject in this short note). Another direction
of arithmetic investigations are the values of~\thetag{1}
at $z=1$ (or $z=-1$), so called {\it zeta values}.
This goes back to Euler's time, who has contributed
by the formula
$$
\Li_{2k}(1)=\zeta(2k)=\frac{(-1)^{k-1}(2\pi)^{2k}B_{2k}}{2(2k)!}
\qquad\text{for}\quad k=1,2,\dots,
\tag2
$$
where $B_{2k}\in\Bbb Q$ are the Bernoulli numbers,
thus Lindemann's proof of the transcendence of~$\pi$
results in the transcendence of the numbers~\thetag{2}.
Ap\'ery~\cite{Ap} (1978) has shown that $\zeta(3)$ is irrational,
and since that time, thanks to Ball and Rivoal~\cite{BR},
we dispose of only partial
irrationality information for other values of
$\Li_{2k+1}(1)=\zeta(2k+1)$ if $k=2,3,4,\dots$\,.

All known achievements in this subject are closely
related to hypergeometric series and also multiple
and complex integrals originated from the series.
This is a general concept of the hypergeometric
method developed for arithmetic study of the
values of the polylogarithms; we refer the reader
to a brief exposition of this concept in~\cite{Ne}.

Here we would like to present some new ingredients
of the hypergeometric method. We cannot achieve some
new number-theoretic results by these means,
and for the moment this note may be viewed as a methodological
contribution. Nevertheless, approximations to the values
of the polylogarithms that we derive here are
quite reasonable from the computational point of view,
and, in this sense, we continue our previous work
on deducing curious Ap\'ery-like recurrences.

We hope that the reader is somehow familiar with
our work on the hypergeometric method in arithmetic
study of zeta values (at least with the preprint~\cite{Zu1}).

\head
1. Simultaneous approximations to the logarithm and dilogarithm
\endhead

For each $n=0,1,\dots$, consider the rational function
$$
R_n(t)=\frac{((t-1)(t-2)\dotsb(t-n))^2}
{n!\cdot t(t+1)\dotsb(t+n)}.
$$
Since degree of the numerator is greater than
degree of the denominator, we will have a polynomial
part while decomposing into partial fractions.
The arithmetic properties of this decomposition
are given in the following statement; $D_n$~denotes
the least common multiple of the numbers $1,2,\dots,n$.

\proclaim{Lemma 1}
We have
$$
R_n(t)=\frac{((t-1)(t-2)\dotsb(t-n))^2}
{n!\cdot t(t+1)\dotsb(t+n)}
=\sum_{k=0}^n\frac{A_k}{t+k}+B(t),
$$
where numbers $A_k$ are all integers and
$D_n\cdot B(t)$ is an integer-valued polynomial
of degree~$n-1$.
\endproclaim

\demo{Proof}
Write this decomposition as follows:
$$
R_n(t)=\sum_{k=0}^n\frac{A_k}{t+k}
+\sum_{j=0}^{n-1}B_j\frac{t(t+1)\dotsb(t+j-1)}{j!}
$$
(the empty product for $j=0$ is~1). The coefficients
$A_k$ are easily determined by the standard procedure:
$$
A_k=R_n(t)(t+k)\big|_{t=-k}
=(-1)^k\binom nk\binom{n+k}k^2,
\qquad k=0,1,\dots,n,
\tag3
$$
while the remaining group of unknown coefficients requires
some work. Denote
$$
F_l(t)=(t+l)\cdot\sum_{k=0}^n\frac{A_k}{t+k},
\qquad l=0,1,\dots,n.
$$
Then
$$
F_l(t)=\sum_{k=0}^nA_k\biggl(1-\frac{k-l}{t+k}\biggr)
$$
and
$$
\frac{\d F_l(t)}{\d t}
=\sum\Sb k=0\\k\ne l\endSb^nA_k\frac{k-l}{(t+k)^2}.
$$
Therefore,
$$
\frac{\d F_l(t)}{\d t}\bigg|_{t=-l}
=\sum\Sb k=0\\k\ne l\endSb^n\frac{A_k}{k-l},
\qquad l=0,1,\dots,n.
\tag4
$$
What will happen if we do the same with the polynomial tail?
Define
$$
G_l(t)=(t+l)\cdot\sum_{j=0}^{n-1}B_j
\frac{t(t+1)\dotsb(t+j-1)}{j!},
\qquad l=0,1,\dots,n-1.
$$
Then
$$
\frac{\d G_l(t)}{\d t}\bigg|_{t=-l}
=\frac{\d}{\d t}\biggl(\sum_{j=0}^lB_j
\frac{t(t+1)\dotsb(t+j-1)}{j!}\cdot(t+l)\biggr)
\bigg|_{t=-l}
$$
and, since for a polynomial~$P(t)$
$$
\frac{\d}{\d t}\bigl(P(t)(t+l)\bigr)\bigg|_{t=-l}
=P(-l),
$$
we finally obtain
$$
\frac{\d G_l(t)}{\d t}\bigg|_{t=-l}
=\sum_{j=0}^l(-1)^j\binom ljB_j,
\qquad l=0,1,\dots,n-1,
$$
and hence
$$
B_l=\frac{\d G_l(t)}{\d t}\bigg|_{t=-l}
-(-1)^l\sum_{j=0}^{l-1}(-1)^j\binom ljB_j,
\qquad l=0,1,\dots,n-1.
\tag5
$$
Furthermore, we have
$$
\gather
F_l(t)+G_l(t)=R_n(t)(t+l)
=\biggl(\frac{(t-1)\dotsb(t-n)}{n!}\biggr)^2
\cdot\frac{n!(t+l)}{t(t+1)\dotsb(t+n)},
\\
l=0,1,\dots,n-1,
\endgather
$$
hence
$$
D_n\cdot\frac{\d}{\d t}\bigl(F_l(t)+G_l(t)\bigr)\bigg|_{t=-l}
\in\Bbb Z,
\qquad l=0,1,\dots,n-1,
$$
where $D_n$ denotes the least common multiple of the numbers
$1,2,\dots,n$. These inclusions and inclusions
$$
D_n\cdot\frac{\d F_l(t)}{\d t}\bigg|_{t=-l}
\in\Bbb Z,
\qquad l=0,1,\dots,n-1,
$$
which follow from formulae \thetag{3} and \thetag{4},
together with the induction on~$l$ on the basis of~\thetag{5},
show that
$$
D_n\cdot B_l\in\Bbb Z,
\qquad l=0,1,\dots,n-1,
$$
and the proof follows.
\enddemo

Since for an integer-valued polynomial $P(t)$
of degree at most~$n$ its derivative multiplied by~$D_n$
is again an integer-valued polynomial, we also have

\proclaim{Lemma 2}
The following decomposition is valid:
$$
-\frac{\d R_n(t)}{\d t}
=\sum_{k=0}^n\frac{A_k}{(t+k)^2}+\wt B(t),
$$
where numbers $A_k$ are all integers and
$D_n^2\cdot\wt B(t)$ is an integer-valued polynomial
of degree~$n-2$.
\endproclaim

Let $z$ be a rational number with $0<|z|<1$.
We are now interested in the following two
hypergeometric-type series:
$$
r_n=r_n(z)
=\sum_{\nu=1}^\infty z^\nu R_n(t)\big|_{t=\nu},
\qquad
\wt r_n=\wt r_n(z)
=-\sum_{\nu=1}^\infty z^\nu\frac{\d R_n(t)}{\d t}\bigg|_{t=\nu}.
$$

\proclaim{Lemma 3}
We have
$$
\gathered
r_n(z)=a_n\Li_1(z)-b_n,
\qquad
\wt r_n(z)=a_n\Li_2(z)-\wt b_n,
\\
z_1^na_n\in\Bbb Z,
\quad
(z_1z_2)^nD_nb_n\in\Bbb Z,
\quad
(z_1z_2)^nD_n^2\wt b_n\in\Bbb Z,
\endgathered
\tag6
$$
where $z_1$ and $z_2$ are the denominators
of the numbers $1/z$ and $z/(1-z)$, respectively.
\endproclaim

\demo{Proof}
Let us write the polynomials $B(t)$ and $\wt B(t)$ in the form
$$
B(t)=\sum_{j=0}^{n-1}B_j\frac{(t-1)(t-2)\dotsb(t-j)}{j!},
\qquad
\wt B(t)=\sum_{j=0}^{n-2}\wt B_j\frac{(t-1)(t-2)\dotsb(t-j)}{j!},
$$
where
$$
D_n\cdot B_j\in\Bbb Z, \quad
D_n^2\cdot\wt B_j\in\Bbb Z,
\qquad j=0,\dots,n-1
\tag7
$$
(this is guaranteed by the theorem of choosing a basis
in the $\Bbb Z$-space of integer-valued polynomials).
Then
$$
\align
r_n
&=\sum_{\nu=1}^\infty z^\nu
\biggl(\sum_{k=0}^n\frac{A_k}{\nu+k}
+\sum_{j=0}^{n-1}B_j
\frac{(\nu-1)(\nu-2)\dotsb(\nu-j)}{j!}\biggr)
\\
&=\sum_{k=0}^nA_kz^{-k}
\sum_{\nu=1}^\infty\frac{z^{\nu+k}}{\nu+k}
+\sum_{j=0}^{n-1}B_jz^{j+1}\sum_{\nu=1}^\infty
\frac{(\nu-1)(\nu-2)\dotsb(\nu-j)}{j!}z^{\nu-j-1}
\allowdisplaybreak
&=\sum_{k=0}^nA_kz^{-k}
\biggl(\sum_{l=1}^\infty-\sum_{l=1}^k\biggr)\frac{z^l}l
-\sum_{j=0}^{n-1}B_jz^{j+1}\cdot\frac1{(z-1)^{j+1}}
\\
&=\sum_{k=0}^nA_kz^{-k}\cdot\Li_1(z)
-\sum_{k=0}^nA_k\sum_{l=1}^k\frac{z^{-(k-l)}}l
-\sum_{j=0}^{n-1}B_j\biggl(\frac z{z-1}\biggr)^{j+1}.
\endalign
$$
In the same vein,
$$
\wt r_n
=\sum_{k=0}^nA_kz^{-k}\cdot\Li_2(z)
-\sum_{k=0}^nA_k\sum_{l=1}^k\frac{z^{-(k-l)}}{l^2}
-\sum_{j=0}^{n-2}\wt B_j\biggl(\frac z{z-1}\biggr)^{j+1}.
$$
Using \thetag{7} and integrality of $A_k$ for
all $k=0,1,\dots,n$, we arrive at the desired claim.
\enddemo

\remark{Remark}
As follows from~\thetag{3} and the above proof,
we have the following explicit formula:
$$
a_n=\sum_{k=0}^n\binom nk\binom{n+k}k^2
\biggl(-\frac1z\biggr)^k.
$$
\endremark

As we see from Lemma~3, the sequences $r_n(z)$ and $\wt r_n(z)$,
$n=0,1,\dots$, realize simultaneous rational approximations
to $\Li_1(z)$ and $\Li_2(z)$, although approximation `tails' are
not simply polynomials in~$1/z$, but sums of two polynomials from
$\Bbb Q[1/z]$ and $\Bbb Q[z/(z-1)]$. Although we worked in the
region $|z|<1$, the result of Lemma~3 remains valid in the closed
disc $|z|\le1$ except the point $z=1$ by analytic continuation.


Running the Gosper--Zeilberger algorithm of creative
telescoping \cite{PWZ, Chapter~6}
with the input $R_n(t)z^t$, one can find
the difference operator of order~3,
which annihilates linear forms~\thetag{6}
and their coefficients $a_n,b_n,\wt b_n$.
In order not to frighten the reader, we indicate
only the characteristic polynomial of this operator
$$
z(z-1)\lambda^3-(3z^2-20z+16)\lambda^2+z(3z+8)\lambda-z^2
$$
(containing all details on the asymptotic behaviour
of the approximants, due to Poincar\'e's theorem),
and the partial case $z=-1$ of the corresponding recurrence.

\proclaim{Theorem 1}
The coefficients $a_n$, $b_n$ and $\wt b_n$ of
the simultaneous approximations
$$
\aligned
r_n
&=r_n(-1)=a_n\Li_1(-1)-b_n=-a_n\log2-b_n,
\\
\wt r_n
&=\wt r_n(-1)=a_n\Li_2(-1)-\wt b_n
=-a_n\frac{\pi^2}{12}-\wt b_n,
\endaligned
\qquad n=0,1,\dots,
$$
as well as the approximations themselves satisfy
the recurrence relation
$$
\gather
\aligned
&
2(59n-24)(n+1)^2a_{n+1}
-(2301n^3+1365n^2-376n-240)a_n
\\ \vspace{-1.5pt} &\qquad
-(295n^3-120n^2-60n+35)a_{n-1}
-(59n+35)(n-1)^2a_{n-2}=0,
\endaligned
\\
n=2,3,\dots,
\endgather
$$
of order~$3$, and the necessary initial data is as follows:
$$
\gathered
a_0=1, \quad a_1=5, \quad a_2=55,
\\
b_0=0, \quad b_1=-\frac72, \quad b_2=-\frac{305}8,
\qquad
\wt b_0=0, \quad \wt b_1=-4, \quad \wt b_2=-\frac{181}4.
\endgathered
$$
In addition,
$$
\gather
\lim_{n\to\infty}|r_n|^{1/n}
=\lim_{n\to\infty}|\wt r_n|^{1/n}
=|\lambda_{1,2}|=0.15960248\dots,
\\
\lim_{n\to\infty}|a_n|^{1/n}
=\lim_{n\to\infty}|b_n|^{1/n}
=\lim_{n\to\infty}|\wt b_n|^{1/n}
=\lambda_3=19.62866250\dots,
\endgather
$$
where $\lambda_{1,2}=-0.06433125\hdots\pm i0.14606314\dots$
and $\lambda_3$ are
zeros of the characteristic polynomial
$2\lambda^3-39\lambda^2-5\lambda-1$.
\endproclaim

One can run the algorithm of creative telescoping with
the input $R_n(t)$ (that is, $z=1$, a non-sense!) to
obtain a much simpler recurrence
$$
(n+1)^2a_{n+1}+(11n^2+11n+3)a_n-n^2a_{n-1}=0,
\qquad n=1,2,\dots,
$$
of order~$2$, which may be recognized in view
of Ap\'ery's proof of the irrationality of~$\zeta(2)$
(see \cite{Ap} and~\cite{Po}).
In fact, we have the identity
$$
a_n=\sum_{k=0}^n(-1)^k\binom nk\binom{n+k}k^2
=(-1)^n\sum_{k=0}^n\binom{n+k}k\binom nk^2
$$
thanks to Thomae's transformation of
${}_3F_2(1)$-hypergeometric series (see \cite{Ba, Section~3.2}),
and the latter sum gives (up to the sign factor)
the denominators of Ap\'ery's approximations to~$\zeta(2)$.
In order to give the necessary sense to the substitution
$z=1$, we should introduce the following complex Barnes integral:
$$
\frac1{2\pi i}\int_{C-i\infty}^{C+i\infty}
R_n(t)\biggl(\frac\pi{\sin\pi t}\biggr)^2z^t\,\d t
=\wt r_n(z)-r_n(z)\log z,
\tag8
$$
where $C$ is an arbitrary constant in the interval
$0<C<n+1$. The integral converges in the whole disc
$|z|\le1$, hence we should have the limit in the
right-hand side of~\thetag{8}, i.e.\ (in notations
of the proof of Lemma~3) the limit
$$
\lim_{z\to1}\biggl(\sum_{j=0}^{n-2}\wt B_j
\biggl(\frac z{z-1}\biggr)^{j+1}
-\log z\cdot\sum_{j=0}^{n-1}B_j
\biggl(\frac z{z-1}\biggr)^{j+1}\biggr)
$$
exists and is equal to a certain rational constant
depending on~$n$.
On the other hand, the complex integral in~\thetag{8}
admits a real double-integral representation thanks
to~\cite{Ne, Theorem~2}. Taking $m=3$, $r=2$ and
$a_1=a_2=a_3=n+1$, $b_2=b_3=2n+2$ in this Nesterenko's
theorem, we obtain
$$
\wt r_n(z)-r_n(z)\log z
=\iint\limits_{[0,1]^2}
\frac{x^n(1-x)^ny^n(1-y)^n}{(1-x+zxy)^{n+1}}\,\d x\,\d y.
$$
Substituting $z=1$ and multiplying by~$(-1)^n$ reduces
the latter integral to Beukers' famous double integral~\cite{Be}
for Ap\'ery's approximations to~$\zeta(2)$.

\head
2. Simultaneous approximations to $\zeta(2)$ and $\zeta(3)$
\endhead

Our first natural generalization of the construction
in the previous section is based on the rational function
$$
R_n(t)=\frac{((t-1)(t-2)\dotsb(t-n))^3}
{n!^2\cdot t(t+1)\dotsb(t+n)}.
$$
Then `reasonable' approximations to the first three
polylogarithms are given by the series
$$
\align
r_n(z)
&=\sum_{\nu=1}^\infty z^\nu R_n(t)\big|_{t=\nu}
=a_n\Li_1(z)-b_n,
\\
\wt r_n(z)
&=-\sum_{\nu=1}^\infty z^\nu
\frac{\d R_n(t)}{\d t}\bigg|_{t=\nu}
=a_n\Li_2(z)-\wt b_n,
\\
\wt{\wt r}_n(z)
&=\frac12\sum_{\nu=1}^\infty z^\nu
\frac{\d^2R_n(t)}{\d t^2}\bigg|_{t=\nu}
=a_n\Li_3(z)-\wt{\wt b}_n,
\endalign
$$
where
$$
z_1^na_n\in\Bbb Z,
\quad
(z_1z_2)^nD_nb_n\in\Bbb Z,
\quad
(z_1z_2)^nD_nD_{2n}\wt b_n\in\Bbb Z,
\quad
(z_1z_2)^nD_nD_{2n}^2\wt{\wt b}_n\in\Bbb Z,
\tag9
$$
and $z_1$ and $z_2$ are the denominators
of the numbers $1/z$ and $z/(1-z)$, respectively.
The reason of having the multiples $D_{2n}$
in~\thetag{9} is the higher degree $2n-1$
of the polynomial in the decomposition of $R_n(t)$
into partial fractions, and it is required to
derivate it for getting the representation
of $\wt r_n$ and $\wt{\wt r}_n$. The explicit
formula for the coefficient~$a_n$ is as follows:
$$
a_n=(-1)^n\sum_{k=0}^n\binom nk\binom{n+k}k^3\biggl(-\frac1z\biggr)^k.
$$

This time we are interested in the particular
`non-sense' case $z=1$ of the construction.
Without tiring reader's eyes by writing two
complex integrals converging in the disc $|z|\le1$
and guaranteeing the existence of limits for
corresponding series expansions, we just
present the final result for the approximation
sequences
$$
\wt r_n=\wt r_n(1)
=a_n\zeta(2)-\wt b_n,
\quad
\wt{\wt r}_n=\wt{\wt r}_n(1)
=a_n\zeta(3)-\wt{\wt b}_n,
\qquad n=0,1,\dotsc.
$$

\proclaim{Theorem 2}
The above sequences $\wt r_n$, $\wt{\wt r}_n$
as well as the coefficients $a_n$, $\wt b_n$ and $\wt{\wt b}_n$
satisfy the recurrence relation
$$
\align
&
2(946n^2-731n+153)(2n+1)(n+1)^3a_{n+1}
\\ &\qquad
-2(104060n^6+127710n^5+12788n^4-34525n^3-8482n^2+3298n+1071)a_n
\\ &\qquad
+2(3784n^5-1032n^4-1925n^3+853n^2+328n-184)na_{n-1}
\\ &\qquad
-(946n^2+1161n+368)n(n-1)^3a_{n-2}=0,
\qquad n=2,3,\dots,
\endalign
$$
of order~$3$, and the necessary initial data is as follows:
$$
\gathered
a_0=1, \quad a_1=7, \quad a_2=163,
\\
\wt b_0=0, \quad \wt b_1=\frac{23}2,
\quad \wt b_2=\frac{2145}8,
\qquad
\wt{\wt b}_0=0, \quad \wt{\wt b}_1=\frac{17}2,
\quad \wt{\wt b}_2=\frac{3135}{16}.
\endgathered
$$
In addition,
$$
\gather
\limsup_{n\to\infty}|\wt r_n|^{1/n}
=\limsup_{n\to\infty}|\wt{\wt r}_n|^{1/n}
=|\lambda_{1,2}|=0.067442248\dots,
\\
\lim_{n\to\infty}|a_n|^{1/n}
=\lim_{n\to\infty}|\wt b_n|^{1/n}
=\lim_{n\to\infty}|\wt{\wt b}_n|^{1/n}
=\lambda_3=54.96369509\dots,
\endgather
$$
where $\lambda_{1,2}=0.018152450\hdots\pm i0.064953409\dots$ and $\lambda_3$ are
zeros of the characteristic polynomial
$4\lambda^3-220\lambda^2+8\lambda-1$.
\endproclaim

Based on the recurrence, we have observed experimentally and
we are able to show that
the correct form of the inclusions~\thetag{9} in this special
case $z=1$ is
$$
\gathered
a_n=(-1)^n\sum_{k=0}^n(-1)^k\binom nk\binom{n+k}k^3
=\sum_{k=0}^n\binom nk^2\binom{n+k}n\binom{n+2k}n\in\Bbb Z,
\\
D_nD_{2n}\wt b_n\in\Bbb Z,
\quad
D_n^3\wt{\wt b}_n\in\Bbb Z.
\endgathered
$$

\remark{Remark}
Normalizing the approximations of Theorem~3
by multiplying them by the factor $\binom{2n}n$,
we arrive at the recurrence previously obtained
in~\cite{Zu2, Theorem~4} by means of a certain
implicit construction. Our new explicit consideration
gives an answer to arithmetic observations
posed in~\cite{Zu2}.
\endremark

\head
3. Well-poised approximations
\endhead

The arithmetic study of zeta values was strongly influenced
by {\it well-poised\/} hypergeometric series.
They are a `heart' of the proof in~\cite{BR} and of many other similar results,
and we cannot avoid considering a well-poised generalization
of the construction in Section~1.

Take
$$
\align
R_n(t)
&=(-1)^{n+1}\biggl(t+\frac n2\biggr)
\frac{((t-1)\dotsb(t-n)\cdot(t+n+1)\dotsb(t+2n))^2}
{n!\cdot(t(t+1)\dotsb(t+n))^3}
\\
&=(-1)^nR_n(-t-n),
\endalign
$$
which has now a quite complicated partial-fraction
decomposition:
$$
R_n(t)
=\sum_{k=0}^n\biggl(\frac{A_k}{(t+k)^3}
+\frac{A_k'}{(t+k)^2}+\frac{A_k''}{t+k}\biggr)
+B(t),
\tag10
$$
where (repeating arguments of the proof of Lemma~1)
$2A_k$, $2D_nA_k'$, $2D_n^2A_k''$ are integers
for all $k=0,1,\dots,n$, while $2D_n^3B(t)$ is
an integer-valued polynomial. Gathering this
arithmetic knowledge and proceeding as in the
proof of Lemma~3, we deduce that the series
$$
r_n(z)
=\sum_{\nu=1}^\infty z^\nu R_n(t)\big|_{t=\nu},
\qquad
\wt r_n(z)
=-\sum_{\nu=1}^\infty z^\nu
\frac{\d R_n(t)}{\d t}\bigg|_{t=\nu}
$$
are certain linear forms involving certain polylogarithms
(up to $\Li_4(z)$). Not so exciting, but we would like
to deal with the construction at the only one point,
$z=-1$. The well-poised thread (take $-t-n$ in place
of~$t$ in~\thetag{10}) results in equalities
$(-1)^kA_k=-(-1)^{n-k}A_{n-k}$ and
$(-1)^kA_k''=-(-1)^{n-k}A_{n-k}''$, $k=0,1,\dots,n$,
and they are the circumstance, which makes
$r_n(-1)$ and $\wt r_n(-1)$ linear forms
in $\Li_2(-1)=-\pi^2/12$, $1$
and $2\Li_3(-1)=-3\zeta(3)/2$, $1$, respectively, with the same
leading coefficient. We write this final production
as follows:
$$
r_n=r_n(-1)=a_n\frac{\pi^2}{12}-b_n,
\quad
\wt r_n=\wt r_n(-1)=a_n\frac{3\zeta(2)}2-\wt b_n,
\qquad n=0,1,\dots,
$$
where
$$
2D_na_n\in\Bbb Z,
\quad
2^nD_n^3b_n\in\Bbb Z,
\quad
2^nD_n^4b_n\in\Bbb Z,
\qquad n=0,1,\dots;
\tag11
$$
the $n$th powers of $2$ appear since the two is the denominator
of $z/(z-1)$ when $z=-1$.
Applying the algorithm of creative telescoping we arrive
at the following result.

\proclaim{Theorem 3}
The sequences $r_n$, $\wt r_n$
and the coefficients $a_n$, $b_n$ and $\wt b_n$
satisfy the recurrence relation
$$
\align
&
(1457n^2-1363n+348)(n+1)^4a_{n+1}
\\ &\qquad
-(148614n^6+158202n^5-9295n^4-61894n^3-11111n^2+8932n+2784)a_n
\\ &\qquad
+(97619n^6-91321n^5-9443n^4+35343n^3-5440n^2-5678n+1768)a_{n-1}
\\ &\qquad
-3(1457n^2+1551n+442)(3n-2)(3n-4)(n-1)^2a_{n-2}=0,
\qquad n=2,3,\dots,
\endalign
$$
of order~$3$, and the initial values are as follows:
$$
\gathered
a_0=1, \quad a_1=8, \quad a_2=264,
\\
b_0=0, \quad b_1=\frac{13}2,
\quad b_2=\frac{1737}8,
\qquad
\wt b_0=0, \quad \wt b_1=\frac{29}2,
\quad \wt b_2=\frac{7617}{16}.
\endgathered
$$
In addition,
$$
\gather
\limsup_{n\to\infty}|r_n|^{1/n}
=\limsup_{n\to\infty}|\wt r_n|^{1/n}
=|\lambda_{1,2}|=0.51616460\dots,
\\
\lim_{n\to\infty}|a_n|^{1/n}
=\lim_{n\to\infty}|b_n|^{1/n}
=\lim_{n\to\infty}|\wt b_n|^{1/n}
=\lambda_3=101.34149804\dots,
\endgather
$$
where $\lambda_{1,2}=0.32925097\hdots\pm i0.39751691\dots$ and $\lambda_3$ are
zeros of the characteristic polynomial
$\lambda^3-102\lambda^2+67\lambda-27$.
\endproclaim

On the basis of the recurrence relation we find
much better inclusions than~\thetag{11}, namely
$$
\wt\Phi_n^{-1}a_n\in\Bbb Z,
\quad
2\wt\Phi_n^{-1}D_n^2b_n\in\Bbb Z,
\quad
2\wt\Phi_n^{-1}D_n^3b_n\in\Bbb Z,
\qquad n=0,1,\dots,
\tag12
$$
where $\wt\Phi_n$ is the following product over primes:
$$
\wt\Phi_n=\prod\Sb p\le n\\2/3\le\{n/p\}<1\endSb p,
$$
$\{\,\cdot\,\}$~denotes the fractional part of a number.
Inclusions~\thetag{12} are quite expected in view
of `denominator conjectures' around linear forms in
zeta values (see \cite{KR, Section~17.1} about the difficulties
in proving the correct arithmetic in similar cases).
But why do we get the cancellation of~$2^n$? This
might be also caused by the well-poised origin of
the series used by us. At least the integrality
of~$a_n$ is an immediate consequence of
the following explicit formulae:
$$
\align
a_n
&=(-1)^n\sum_{j=0}^n\frac{\d}{\d t}\biggl(\frac n2-t\biggr)
\binom nt^3\binom{n+t}n^2\binom{2n-t}n^2\biggr|_{t=j}
\\
&=\mathop{\sum\sum}_{0\le i\le j\le n}
(-1)^{n+j}\binom ni^2\binom nj\binom{2n-i}n\binom{n+j}n
\binom{n+j-i}n
\tag13
\\
&=\sum_{i=0}^n\sum_{j=0}^n
\binom ni^2\binom nj^2\binom{n+i}n\binom{i+j}i,
\tag14
\endalign
$$
where equality~\thetag{13} follows
from~\cite{KR, Proposition~5} and equality~\thetag{14} is
communicated to us by G.~Almkvist.

\head
4. Final remarks
\endhead

As already promised, no new irrationality and linear
independence results were presented. We just have tried
to give some sense to certain hypergeometric-type series that
are expressed in terms of polylogarithms and are
divergent when one formally plugs $z$ with $|z|=1$.
Transforming a non-terminating single hypergeometric series
into a multiple one (some kind of
`identit\'es non-termin\'ees gigantesques',
cf\. \cite{KR, Section~17.5})
often meets convergence troubles
for the latter series, i.e\. it is just a formal
transformation, which we could never use in a rigorous
proof. As an option to proceed in such troubling
cases, we see dealing with transformations for
complex Barnes (multiple) integrals and further
decompositions of the integrals into sums
involving (multiple) zeta values. This does not
look an easy program, but we do not believe that
deducing new results for zeta values and polylogarithms
might be simple.



\Refs
\widestnumber\key{WW}

\ref\key Ap
\by R.~Ap\'ery
\paper Irrationalit\'e de $\zeta(2)$ et $\zeta(3)$
\jour Ast\'erisque
\vol61
\yr1979
\pages11--13
\endref

\ref\key Ba
\by W.\,N.~Bailey
\book Generalized hypergeometric series
\bookinfo Cambridge Math. Tracts
\vol32
\publ Cambridge Univ. Press
\publaddr Cambridge
\yr1935
\moreref
\bookinfo 2nd reprinted edition
\publaddr New York
\publ Stechert-Hafner
\yr1964
\endref

\ref\key BR
\by K.~Ball and T.~Rivoal
\paper Irrationalit\'e d'une infinit\'e de valeurs de la
fonction z\^eta aux entiers impairs
\jour Invent. Math.
\vol146
\issue1
\yr2001
\pages193--207
\endref

\ref\key Be
\by F.~Beukers
\paper A note on the irrationality of~$\zeta(2)$ and~$\zeta(3)$
\jour Bull. London Math. Soc.
\vol11
\issue3
\yr1979
\pages268--272
\endref

\ref\key Ch
\by G.\,V.~Chudnovsky
\paper Pad\'e approximations to the generalized
hypergeometric functions. I
\jour J. Math. Pures Appl. (9)
\vol58
\issue4
\yr1979
\pages445--476
\endref

\ref\key Ga
\by A.\,I.~Galochkin
\paper Lower estimates of linear forms in values of certain $G$-functions
\jour Mat. Zametki [Math. Notes]
\vol18
\issue4
\yr1975
\pages541--552
\endref

\ref\key Ha1
\by M.~Hata
\paper On the linear independence of the values of
polylogarithmic functions
\jour J. Math. Pures Appl. (9)
\vol69
\issue2
\yr1990
\pages133--173
\endref

\ref\key Ha2
\by M.~Hata
\paper Rational approximations to the dilogarithm
\jour Trans. Amer. Math. Soc.
\vol336
\issue1
\yr1993
\pages363--387
\endref

\ref\key KR
\by C.~Krattenthaler and T.~Rivoal
\paper Hyperg\'eom\'etrie et fonction z\^eta de Riemann
\jour Preprint (February 2004)
\finalinfo E-print {\tt math.NT/0311114v3}
\endref

\ref\key Ma
\by W.~Maier
\paper Potenzreihen irrationalen Grenzwertes
\jour J. Reine Angew. Math.
\vol156
\yr1927
\pages93--148
\endref

\ref\key Ne
\by Yu.\,V.~Nesterenko
\paper Integral identities and constructions of approximations
to zeta-values
\jour J. Th\'eorie Nombres Bordeaux
\vol15
\yr2003
\issue2
\pages535--550
\endref

\ref\key Ni
\by E.\,M.~Nikishin
\paper On irrationality of values of functions $F(x,s)$
\jour Mat. Sb. [Russian Acad. Sci. Sb. Math.]
\vol109
\yr1979
\issue3
\pages410--417
\endref

\ref\key PWZ
\by M.~Petkov\v sek, H.\,S.~Wilf and D.~Zeilberger
\book $A=B$
\publaddr Wellesley, MA
\publ A.\,K.~Peters, Ltd.
\yr1996
\endref

\ref\key Po
\by A.~van der Poorten
\paper A proof that Euler missed...
Ap\'ery's proof of the irrationality of~$\zeta(3)$
\paperinfo An informal report
\jour Math. Intelligencer
\vol1
\issue4
\yr1978/79
\pages195--203
\endref

\ref\key RV
\by G.~Rhin and C.~Viola
\paper The permutation group method for the dilogarithm
\inbook Pr\'epublication no.~2004-13
\publ D\'epartement de Math.
\publaddr Metz
\yr2004
\endref

\ref\key Zu1
\by W.~Zudilin
\paper An elementary proof of Ap\'ery's theorem
\jour Preprint (February 2002)
\finalinfo E-print {\tt math.NT/0202159}
\endref

\ref\key Zu2
\by W.~Zudilin
\paper A third-order Ap\'ery-like recursion for $\zeta(5)$
\jour Mat. Zametki [Math. Notes]
\vol72
\issue5
\yr2002
\pages733--737
\moreref
\finalinfo E-print {\tt math.NT/0206178}
\endref

\endRefs
\enddocument